\documentclass[12pt, reqno]{amsart}
\usepackage{amsmath, amstext, amsbsy, amssymb, amscd}

\newtheorem{theorem}{Theorem}
\newtheorem{lemma}[theorem]{Lemma}

\theoremstyle{definition}

\theoremstyle{remark}
\newtheorem{remark}[theorem]{Remark}

\newcommand{\Bqt}{B_\la(q,t)}
\newcommand{\BBqt}{\widehat{B}_\la(q,t)}
\newcommand{\BBo}{\widehat{B}_\emptyset(q,t)}
\newcommand{\C}{ \mathbb C }

\newcommand{\End}{{\rm End}}

\newcommand{\la}{\lambda}

\newcommand{\Q}{ \mathbb Q }
\newcommand{\Tr}{ {\rm Tr} }
\newcommand{\vac}{|0\rangle}

\newcommand{\Z}{ \mathbb Z }

{\vskip-\lastskip\medskip
  \noindent
  {\em #1.}\enspace
  }%
{\qed\par\medskip
  }

\begin{document}
\title[] {The correlation functions of vertex operators and Macdonald polynomials}

\author[Shun-Jen Cheng]{Shun-Jen Cheng}
\address{Department of Mathematics, National
Taiwan University, Taipei, Taiwan
106}\email{chengsj@math.ntu.edu.tw}

\author[Weiqiang Wang]{Weiqiang Wang}
\address{Department of Mathematics, University of Virginia,
Charlottesville, VA 22904} \email{ww9c@virginia.edu}
%\thanks{Partially supported by NSA and NSF}

\maketitle
\date{}
%\tableofcontents

\begin{abstract}
The $n$-point correlation functions introduced by Bloch and
Okounkov have already found several geometric connections and
algebraic generalizations. In this Note we formulate a
$q$,$t$-deformation of this $n$-point function. The key operator
used in our formulation arises from the theory of Macdonald
polynomials and affords a vertex operator interpretation. We
obtain closed formulas for the $n$-point functions when $n =1,2$
in terms of the basic hypergeometric functions. We further
generalize the $q$,$t$-deformed $n$-point function to more general
vertex operators.
\end{abstract}

\section{Introduction}

In \cite{BO} Bloch and Okounkov formulated an $n$-point
correlation function on a Fock space and established a remarkable
closed formula in terms of theta functions (also cf. \cite{Ok}).
Recently, this $n$-point function has found geometric connections
in terms of Gromov-Witten theory \cite{OP} and Hilbert schemes of
points \cite{LQW}, and it also affords several other algebraic
generalizations (cf. \cite{Mi, Wa, CW}). The formulation in
\cite{BO, Ok} boils down to a remarkable operator $T(t)$ on the
ring of symmetric functions which diagonalizes the Schur functions
with explicit eigenvalues.

In this Note we formulate a deformed version of the $n$-point
functions of Bloch-Okounkov, denoted by $\widehat{F}(q_1, t_1;
\ldots; q_n, t_n)$, which also depends on an indeterminate $v$
associated to the energy operator. The role of $T(t)$ is replaced
by an operator $\widehat{\mathfrak B}_{q,t}$ (cf.~Garsia-Haiman
\cite{GH}; see Section~\ref{sec:npoint}) which diagonalizes the
modified Macdonald polynomials $\tilde{H}_\la(q,t)$ and affords a
vertex operator interpretation. In Section~\ref{sec:formulas} we
compute the $1$-point function as
 \begin{align*}
  \widehat{F}(q,t)
 =  \frac{(vqt)_\infty}{(t )_\infty(q )_\infty}
 \end{align*}
where $(a)_\infty :=\prod_{i=0}^\infty (1-a v^i)$. We further
found closed formulas for the $2$-point functions in terms of
basic hypergeometric series (Theorem~\ref{th:2point}).

From the viewpoint of vertex operators, it is also possible to
further generalize the notion of the $n$-point function above, and
we compute explicitly some cases in Section~\ref{sec:general}
(Theorems~\ref{th:1pointgeneral} and \ref{th:npoint}). We end this
Note in Section~\ref{sec:discuss} with a discussion of open
problems and possible connections.

\section{Formulation of the $n$-point functions}
\label{sec:npoint}

\subsection{The operators $\mathfrak B_{q,t}$ and $\widehat{\mathfrak B}_{q,t}$}

Let $t, q$ be two indeterminates. Given a partition $\lambda
=(\la_1, \la_2, \ldots)$ of $n$, we denote by $a'(\square)$
 and $l'(\square)$ the {\em coarm} and {\em coleg} of a
given cell $\square$ \cite{Mac, GH}, and denote

$$\Bqt :=\sum_{\square \in \la} q^{a'(\square) } t^{l'(\square) }.$$

We set
\begin{eqnarray*}
%\Wqt &:=& \sum_{i \ge 1} t^{i-1} q^{\la_i - 1/2}  \\
\BBqt &:=& \frac{1}{1-q}\sum_{i \ge 1} t^{i-1} q^{\la_i}.
\end{eqnarray*}
\begin{lemma} \label{lem:convert}
We have
\begin{eqnarray*}
\Bqt
%= \frac{\Wqt -W_\emptyset (q,t)}{q^{1/2} -q^{-1/2}}
 = \BBo -\BBqt
\end{eqnarray*}
where $\BBo =\frac1{(1-q)(1-t)}$.
\end{lemma}
\begin{proof}
We calculate that
\begin{eqnarray*}
\Bqt & =& \sum_{i\ge 1} \frac{t^{i-1} (1-q^{\la_i})}{1-q} \\
     & =& \frac1{(1-q)(1-t)} -\sum_{i\ge 1} \frac{t^{i-1} q^{\la_i}}{1-q}\\
%   &=& \frac1{(1-q)(1-t)} +\frac{\Wqt}{ q^{1/2} -q^{-1/2}} \\
% &=& \frac{\Wqt -W_\emptyset (q,t)}{q^{1/2} -q^{-1/2}} \\
  &=& \BBo -\BBqt.
\end{eqnarray*}
\end{proof}
Note that
\begin{eqnarray}  \label{qt-sym}
\Bqt =B_{\la'}(t,q), \qquad \BBqt ={\widehat B}_{\la'}(t,q).
\end{eqnarray}

Denote by $\Lambda_{q,t}$ the ring of symmetric functions with
coefficients in $\Q (q,t)$. Recall the Macdonald symmetric
functions $P_\la (x; q,t)$, $Q_\la (x; q,t)$ from \cite{Mac} and
its normalized form $J_\la (x; q,t), H_\la (x; q,t),$ and
$\widetilde{H}_\la (x; q,t)$ as in \cite[(8)--(11)]{GH}. We define
the linear operators $\mathfrak B_{q,t}$ and $\widehat{\mathfrak
B}_{q,t}$ on $\Lambda_{q,t}$ (compare \cite[(73), 74)]{GH}) by
letting
\begin{eqnarray*}
%\mathfrak W_{q,t} \widetilde{H}_\la (x; q,t) &=& \Wqt
%\widetilde{H}_\la (x; q,t),    \\
%
\mathfrak B_{q,t} \widetilde{H}_\la (x; q,t) &=& \Bqt
\widetilde{H}_\la (x; q,t),   \\
\widehat{\mathfrak B}_{q,t} \widetilde{H}_\la (x; q,t) &=& \BBqt
\widetilde{H}_\la (x; q,t),  \quad \text{for all } \la.
\end{eqnarray*}

\subsection{The definition of $n$-point correlation functions}
Let $v$ be an indeterminate.  For our purposes we can also think
of $v$ as a complex number with $|v|<1$. For $r \ge 1$ we set
$$(a)_0:=1;\quad (a)_r := \prod_{i=0}^{r-1} (1-a v^i); \quad
(a)_\infty :=\prod_{i=0}^\infty (1-a v^i). $$
The {\em energy operator} $L_0$ on $\Lambda_{q,t}$ is the linear
operator such that $L_0 g= n g$ for every $n$ and every symmetric
function $g$ of degree $n$. Given $\mathfrak f \in \End
(\Lambda_{q,t})$, we consider the trace function
$$\Tr_v \mathfrak f  :=\Tr (v^{L_0} \mathfrak f).$$
In particular for the identity map $I$ we have
$$\Tr_v I  = (v)_\infty^{-1}.$$
The {\em $n$-point (correlation) functions} are defined to be
\begin{eqnarray*}
F(q_1, t_1; \ldots; q_n, t_n)
 := \Tr_v (\mathfrak B_{q_1,t_1} \cdots \mathfrak
B_{q_n,t_n}), \\
\widehat{F}(q_1, t_1; \ldots; q_n, t_n)
 := \Tr_v (\widehat{\mathfrak B}_{q_1,t_1} \cdots \widehat{\mathfrak B}_{q_n,t_n}).
\end{eqnarray*}
We can easily convert between $F$ and $\widehat{F}$ by
Lemma~\ref{lem:convert}.

There is yet another viewpoint. Let $\mathcal P$ be the set of all
partitions, and let $f(\la)$ be a function on $\mathcal P$. We
define the {\em $v$-expectation value} of $f$ to be
$$\langle f \rangle_v := (v)_\infty \sum_{\la \in \mathcal P}
f(\la) v^{|\la|},$$
assuming its convergence.

\begin{lemma}  \label{fun=value}
We have
$$F(q_1, t_1; \ldots, q_n, t_n) = (v)_\infty^{-1}\left \langle \prod_{k=1}^n B_\la(q_k,
t_k) \right \rangle_v.$$
The same relation holds with $F$ and $B$ replaced by $\widehat{F}$
and $\widehat{B}$.
\end{lemma}
\begin{proof}
Note that the operators $\mathfrak B_{q_k,t_k}$ for different $k$
do not commute. Let $\{s_\la\}$ be the Schur functions, cf.
\cite{Mac}, and write $s_\mu = \sum_{\la} a^{(i)}_{\la,\mu}
\widetilde{H}_\la (x; q_i,t_i)$ with $[a^{(i)}_{\la,\mu}]$ being a
triangular matrix with respect to the dominance order. Then
$$\mathfrak B_{q_i,t_i} s_\la = B_\la(q_i,t_i) s_\la + \text{ lower
terms},$$
and thus
$$\mathfrak B_{q_1,t_1}
\cdots \mathfrak B_{q_n,t_n} s_\la = B_\la(q_1,t_1) \cdots
B_\la(q_n,t_n) s_\la + \text{ lower terms}.$$
Therefore,
\begin{eqnarray*} \label{eq:value}
\Tr_v (\mathfrak B_{q_1,t_1}   \cdots \mathfrak B_{q_n,t_n}) =
\sum_\la  B_\la(q_1,t_1) \cdots B_\la(q_n,t_n) v^{|\la|}.
\end{eqnarray*}
\end{proof}
Thanks to (\ref{qt-sym}) and Lemma~\ref{fun=value}, we see that
$F$ and $\widehat{F}$ are symmetric with respect to the
hyperoctahedral group $\Z_2^n \rtimes S_n$, where the symmetric
group $S_n$ permutes the indices $i$ in the pairs $(q_i,t_i)$ and
the $i$-th copy of $\Z_2$ permutes $q_i$ and $t_i$.

\begin{remark}
When $t =q^{-1}$, $\widehat{F}$ reduces to (up to a normalization)
the $n$-point functions introduced by Bloch and Okounkov
\cite{BO}, where the interpretation as a $v$-expectation value was
also made.
\end{remark}

\section{The formulas for $n$-point functions}
\label{sec:formulas}
\subsection{The $1$-point function}

\begin{lemma} \cite[Lemma 6.6]{BO} \label{lem:bo}
For a given $i \ge 1$, we have
$$\langle q^{\la_i} \rangle_v = \frac{(v^i )_\infty}{(v^i q )_\infty}. $$
\end{lemma}
\begin{proof}
By conjugation symmetry of partitions and $\la_i' =\#\{k|\la_k \ge
i\}$, we have
$$ \langle  q^{\la_i} \rangle_v
 =\langle  q^{\la_i'} \rangle_v
 =\frac{(v )_\infty}{(v )_{i-1}(v^i q )_\infty}
 = \frac{(v^i )_\infty}{(v^i q )_\infty}. $$
\end{proof}
We will use for several times the so-called $q$-binomial theorem
(cf. \cite[Appendix II.3]{GR}):
$$
 \sum_{r=0}^\infty t^r \frac{(a)_r}{(v )_r}
 = \frac{(at)_\infty}{(t )_\infty}, \qquad   |t|<1.$$
\begin{theorem} \label{th:onepoint}
The 1-point function is given by:
 \begin{eqnarray*}
 \widehat{F}(q,t) =
 \frac{(vqt)_\infty}{(q )_\infty (t )_\infty}.
 \end{eqnarray*}
 \end{theorem}

\begin{proof}
 We calculate by Lemma~\ref{lem:bo} and the $q$-binomial theorem
 that
 \begin{eqnarray*}
 \left \langle \widehat{B}_\la (q,t) \right \rangle_v
 &=& (1-q)^{-1} \sum_{i=1}^\infty   t^{i-1}
 \frac{(v^i )_\infty}{(v^i q )_\infty} \\
 &=& \frac{(v )_\infty}{(q )_\infty}
 \sum_{r=0}^\infty t^r \frac{(vq )_r}{(v )_r}
 = \frac{( v )_\infty (vqt  )_\infty}{(q )_\infty (t )_\infty}.
\end{eqnarray*}
 By Lemma~\ref{fun=value}, this gives rise to the $1$-point
 function $\widehat{F}(q,t).$
% which is valid for $|t|<1$ and $|v|<1$.
\end{proof}

\begin{remark}
When $t =q^{-1}$, Theorem~\ref{th:onepoint} specializes to
\cite[Theorem~6.5]{BO}.
% Note that we can reformulate Theorem~\ref{th:onepoint} as:
%$
%\begin{eqnarray*}
%\left \langle \widehat{B}_\la (q,t) \right \rangle_v
 % %
% &=& \frac{1}{(1-q)(1-t)}
%  \exp \left ( -\sum_{r=1}^\infty \frac{v^r}{r}
%  \frac{(1-q^r)(1-t^r)}{1-v^r}
%       \right ).
% \end{eqnarray*}
\end{remark}
\subsection{The $2$-point function}

We begin with some preparation.
\begin{lemma} \label{lem:bo2}
For fixed $1 \le i<j$, we have
$$\left \langle  q_1^{\la_i} q_2^{\la_j}  \right \rangle_v
 =\frac{(v)_\infty}{(v)_{i-1} (v^iq_1)_{j-i} (v^j q_1q_2)_\infty}
 =
 \frac{(v)_\infty}{(vq_1q_2)_\infty}
 \frac{(vq_1)_{i-1} (vq_1 q_2)_{j-1}}{ (v)_{i-1}(vq_1)_{j-1}}.
$$
\end{lemma}
\begin{proof}
This is a variant of a special case of \cite[(7.1)]{BO}. Similar
to Lemma~\ref{lem:bo}, it follows directly from $\left \langle
q_1^{\la_i} q_2^{\la_j} \right \rangle_v =\left \langle
q_1^{\la_i'} q_2^{\la_j'}  \right \rangle_v$.
\end{proof}
Set
\begin{eqnarray*}
T_1 &:=&  \frac{(v)_\infty}{(vq_1q_2)_\infty}
 \sum_{i=0}^\infty t_1^i \frac{(vq_1)_i}{(v)_i}
 \sum_{j=i+1}^\infty t_2^j \frac{(vq_1q_2)_j}{(vq_1)_j} \\
 T_2 &:=&   \frac{(v)_\infty}{(vq_1q_2)_\infty}
 \sum_{i=0}^\infty t_2^i \frac{(vq_2)_i}{(v)_i}
  \sum_{j=i+1}^\infty t_1^j \frac{(vq_1q_2)_j}{(vq_2)_j}  \\
T_3 &:=&   \frac{(v)_\infty}{(vq_1q_2)_\infty}
 \sum_{i=0}^\infty (t_1t_2)^i \frac{(vq_1q_2)_i}{(v)_i}.
\end{eqnarray*}

\begin{lemma} \label{threeparts}
We have %
$$\left \langle  \widehat{B}_\la(q_1, t_1)
\widehat{B}_\la(q_2, t_2) \right \rangle_v =(1-q_1)^{-1}
(1-q_2)^{-1} (T_1 +T_2 +T_3).$$
\end{lemma}
\begin{proof}
By definition, we have
\begin{eqnarray*}
&&\left \langle  \widehat{B}_\la(q_1, t_1) \widehat{B}_\la(q_2,
t_2) \right \rangle_v \\
 &=& (1-q_1)^{-1} (1-q_2)^{-1} \left \langle  \sum_{i,j=1}^\infty
t_1^{i-1}q_1^{\la_i} t_2^{j-1}q_2^{\la_j}\right \rangle_v \\
 &=& (1-q_1)^{-1} (1-q_2)^{-1}
 \sum_{i,j=1}^\infty t_1^{i-1} t_2^{j-1}\left \langle
 q_1^{\la_i} q_2^{\la_j}\right \rangle_v \\
 &=& (1-q_1)^{-1} (1-q_2)^{-1}
 \left(\sum_{i<j} + \sum_{i>j} + \sum_{i=j} \right) t_1^{i-1} t_2^{j-1}\left \langle
 q_1^{\la_i} q_2^{\la_j}\right \rangle_v
%
% &=& (1-q_1)^{-1} (1-q_2)^{-1} (T_1 +T_2 +T_3)
\end{eqnarray*}
where the last three summands can be further identified with $T_1,
T_2$ and $T_3$, respectively, using Lemmas~\ref{lem:bo} and
\ref{lem:bo2}.
\end{proof}

For $r\ge 0$, $a_1,\cdots,a_{r+1}\in\C$ and $b_1,\cdots,b_r\in\C$
the $(r+1,r)$-basic hypergeometric series  is the series:
\begin{align*}
{}_{r+1} \Phi_r \left (
\begin{array}{ccc}
  a_1, \cdots, a_{r+1} \\
  b_1,\cdots,b_r
 \end{array}; v; z
 \right):=\sum_{m\ge 0}\frac{(a_1)_m (a_2)_m\cdots(a_{r+1})_m}{(v)_m(b_1)_m\cdots
 (b_r)_m} z^m.
\end{align*}
It is assumed that the denominator is never zero, in which case it
is known to converge absolutely for $|z|<1$ (cf.~\cite{GR}).

\begin{theorem} \label{th:2point}
The $2$-point function $\widehat{F}(q_1, t_1;  q_2, t_2)$ is equal
to
\begin{align*}
&\frac1{(1-q_1) (1-q_2) (1 -t_1t_2)} \cdot
 \frac{
 (vq_1q_2t_1t_2)_\infty}{(vt_1t_2)_\infty(vq_1q_2)_\infty} \cdot
 \\
 &\cdot \left[
 \frac{q_1q_2t_1t_2 -1}{(1-q_1t_1)(1-q_2t_2)}
 + \frac1{1-q_1t_1}
  \;{}_3 \Phi_2 \left ( \begin{array}{ccc}
  v, q_1t_1, vq_1q_2 \\
  vq_1, vq_1q_2t_1t_2
 \end{array}; v; t_2
 \right) \right. \\
& \qquad + \left. \frac1{1-q_2t_2}
  \;{}_3 \Phi_2 \left ( \begin{array}{ccc}
  v, q_2t_2, vq_1q_2 \\
  vq_2, vq_1q_2t_1t_2
 \end{array}; v; t_1
 \right) \right].
\end{align*}
\end{theorem}

\begin{theorem} \label{th:2special}
If $q_1q_2t_1t_2 =1$, then the $2$-point function
$\widehat{F}(q_1, t_1; q_2, t_2)$ is equal to
\begin{align*}
 & \frac1{(1-q_1) (1-q_2) (1 -t_1t_2)} \cdot
 \frac{(v)_\infty}{(vt_1t_2)_\infty (vq_1q_2)_\infty}\cdot \\
 & \cdot \left[ \frac1{1 -q_1t_1}
 \frac{(vt_1^{-1})_\infty (q_2^{-1})_\infty}{(vq_1)_\infty
 (t_2)_\infty}
 +  \frac1{1 -q_2t_2}
 \frac{(vt_2^{-1})_\infty (q_1^{-1})_\infty}{(vq_2)_\infty
 (t_1)_\infty}
 \right]
\end{align*}
\end{theorem}

\begin{proof}[Proof of Theorems~\ref{th:2point} and
\ref{th:2special}] To compute the $2$-point function it suffices
to compute the $T_i$ by Lemma~\ref{threeparts}. First of all,
\begin{eqnarray*}
T_1
 &=&
 %(t_1t_2)^{1/2}
 \frac{(v)_\infty}{(vq_1)_\infty}
 \sum_{i=0}^\infty t_1^i \frac{(vq_1)_i}{(v)_i}
 \sum_{j=i+1}^\infty t_2^j \frac{(v^{j+1}q_1)_\infty}{(v^{j+1}q_1q_2)_\infty} \\
 &=& \frac{(v)_\infty}{(vq_1)_\infty}
 \sum_{i=0}^\infty t_1^i \frac{(vq_1)_i}{(v)_i}
 \sum_{j=i+1}^\infty t_2^j
 \sum_{m=0}^\infty \frac{(q_2^{-1})_m}{(v)_m} (v^{j+1}q_1q_2)^m  \\
 &=& \frac{(v)_\infty}{(vq_1)_\infty}
 \sum_{i=0}^\infty t_1^i \frac{(vq_1)_i}{(v)_i}
 \sum_{m=0}^\infty \frac{(q_2^{-1})_m}{(v)_m} (vq_1q_2)^m
 \sum_{j=i+1}^\infty t_2^j v^{jm} \\
 &=& \frac{(v)_\infty}{(vq_1)_\infty}
 \sum_{i=0}^\infty t_1^i \frac{(vq_1)_i}{(v)_i}
 \sum_{m=0}^\infty \frac{(q_2^{-1})_m}{(v)_m} (vq_1q_2)^m
 \frac{(t_2 v^{m})^{i+1}}{1-t_2v^m} \\
 &=& \frac{(v)_\infty}{(vq_1)_\infty}
 \sum_{m=0}^\infty \frac{(q_2^{-1})_m}{(v)_m} (vq_1q_2)^m
 \frac{t_2 v^{m}}{1-t_2v^m}
 \sum_{i=0}^\infty (t_1t_2 v^{m})^i \frac{(vq_1)_i}{(v)_i} \\
 &=& \frac{(v)_\infty}{(vq_1)_\infty}
 \sum_{m=0}^\infty \frac{(q_2^{-1})_m}{(v)_m} (vq_1q_2)^m
 \frac{t_2 v^{m}}{1-t_2v^m}
 \frac{(v^{m+1}q_1t_1t_2)_\infty}{(v^mt_1t_2)_\infty} \\
 &=& \frac{(v)_\infty(vq_1t_1t_2)_\infty}{(vq_1)_\infty(t_1t_2)_\infty}
 \sum_{m=0}^\infty \frac{(q_2^{-1})_m(t_1t_2)_m}{(v)_m (vq_1t_1t_2)_m}
 (v^2q_1q_2)^m \frac{t_2}{1-t_2v^m}  \\
 &=& \frac{t_2}{1-t_2} \frac{(v)_\infty(vq_1t_1t_2)_\infty}{(vq_1)_\infty(t_1t_2)_\infty}
 \sum_{m=0}^\infty \frac{(q_2^{-1})_m(t_1t_2)_m (t_2)_m}{(v)_m (vq_1t_1t_2)_m (vt_2)_m}
 (v^2q_1q_2)^m.
 \end{eqnarray*}
Thus we obtain that
\begin{eqnarray} \label{hyper}
T_1=\frac{t_2}{1-t_2}
\frac{(v)_\infty(vq_1t_1t_2)_\infty}{(vq_1)_\infty(t_1t_2)_\infty}
 \cdot \; {}_3 \Phi_2 \left (
 \begin{array}{ccc}
 t_1t_2, t_2, q_2^{-1} \\
 vt_2, vq_1t_1t_2
 \end{array}; v; v^2q_1q_2
 \right).
\end{eqnarray}
Here ${}_3 \Phi_2 \left (
 \begin{array}{ccc}
 t_1t_2, t_2, q_2^{-1} \\
 vt_2, vq_1t_1t_2
 \end{array}; v; v^2q_1q_2
 \right)$ is a $(3,2)$-hypergeometric series of type II, since
$$\frac{(vt_2)(vq_1t_1t_2)}{(t_1t_2)(t_2)(q_2^{-1})} =v^2q_1q_2.$$
Recall Hall's transformation formula (\cite{GR}, Appendix III.10):
\begin{align*}
 {}_3 \Phi_2 \left (
 \begin{array}{ccc}
 a, b, c \\
 d, e
 \end{array}; v; \frac{de}{abc}
 \right)
 = \frac{(b)_\infty (de/ab)_\infty (de/bc)_\infty}{(d)_\infty
 (e)_\infty (de/abc)_\infty}
 %\begin{align*}
 {}_3 \Phi_2 \left (
 \begin{array}{ccc}
 d/b, e/b, de/abc \\
 de/ab, de/bc
 \end{array}; v; b
 \right).
\end{align*}
The above two-term transformation formula holds for $|b|<1$ and
$|de/abc|<1$. Applying this transformation formula to
(\ref{hyper}) and cancelling terms with $(t_2)_\infty =(1-t_2)
(vt_2)_\infty$, we rewrite (\ref{hyper}) as
\begin{eqnarray*}
 T_1 &=& t_2
 \frac{(v)_\infty(v^2 q_1)_\infty (v^2q_1q_2t_1t_2)_\infty}{(vq_1)_\infty(t_1t_2)_\infty(v^2q_1q_2)_\infty}
  {}_3 \Phi_2 \left ( \begin{array}{ccc}
 v, vq_1t_1, v^2q_1q_2 \\
 v^2q_1, v^2q_1q_2t_1t_2
 \end{array}; v; t_2
 \right)  \\
 &=& \frac{t_2 }{1-vq_1}
 \frac{(v)_\infty (v^2
 q_1q_2t_1t_2)_\infty}{(t_1t_2)_\infty(v^2q_1q_2)_\infty} \cdot
  {}_3 \Phi_2 \left ( \begin{array}{ccc}
 v, vq_1t_1, v^2q_1q_2 \\
 v^2q_1, v^2q_1q_2t_1t_2
 \end{array}; v; t_2
 \right)  \\
 &=& \frac{t_2 }{1-vq_1}
 \frac{1-vq_1q_2}{1-vq_1q_2t_1t_2}
 \frac{(v)_\infty (vq_1q_2t_1t_2)_\infty}{(t_1t_2)_\infty(vq_1q_2)_\infty} \cdot
  {}_3 \Phi_2 \left ( \begin{array}{ccc}
 v, vq_1t_1, v^2q_1q_2 \\
 v^2q_1, v^2q_1q_2t_1t_2
 \end{array}; v; t_2
 \right)  \\
 &=& \frac{1}{1-q_1t_1}
 \frac{(v)_\infty (vq_1q_2t_1t_2)_\infty}{(t_1t_2)_\infty(vq_1q_2)_\infty} \cdot
 \sum_{m=0}^\infty \frac{(v)_{m+1} (q_1t_1)_{m+1}
 (vq_1q_2)_{m+1}}{(v)_{m+1} (vq_1)_{m+1} (vq_1q_2t_1t_2)_{m+1}}
 t_2^{m+1}  \\
 &=& \frac{1}{1-q_1t_1}
 \frac{(v)_\infty (vq_1q_2t_1t_2)_\infty}{(t_1t_2)_\infty(vq_1q_2)_\infty} \cdot
  \left[
  {}_3 \Phi_2 \left ( \begin{array}{ccc}
  v, q_1t_1, vq_1q_2 \\
  vq_1, vq_1q_2t_1t_2
 \end{array}; v; t_2
 \right) -1 \right] \\
 &=& \frac1{1 -t_1t_2}
 \frac{(v)_\infty (vq_1q_2t_1t_2)_\infty}{(vt_1t_2)_\infty(vq_1q_2)_\infty}
 \cdot\\
 &&\qquad \quad \cdot \left[ \frac1{1-q_1t_1}
  \;{}_3 \Phi_2 \left ( \begin{array}{ccc}
  v, q_1t_1, vq_1q_2 \\
  vq_1, vq_1q_2t_1t_2
 \end{array}; v; t_2
 \right) - \frac1{1-q_1t_1} \right].
\end{eqnarray*}

Since $T_2$ is the same as $T_1$ after switching of variables
$t_1\leftrightarrow t_2, q_1\leftrightarrow q_2$, we have
\begin{eqnarray*}
T_2
  &=&
  \frac1{1 -t_1t_2} \cdot
 \frac{(v)_\infty
 (vq_1q_2t_1t_2)_\infty}{(vt_1t_2)_\infty(vq_1q_2)_\infty} \cdot
 \\
&& \qquad \quad \cdot
 \left [ \frac1{1-q_2t_2}
  \;{}_3 \Phi_2 \left ( \begin{array}{ccc}
  v, q_2t_2, vq_1q_2 \\
  vq_2, vq_1q_2t_1t_2
 \end{array}; v; t_1
 \right) - \frac1{1-q_2t_2} \right].
\end{eqnarray*}

Note in addition by the $q$-binomial theorem that
\begin{eqnarray*}
T_3 = \frac1{1 -t_1t_2}\cdot
 \frac{(v)_\infty
 (vq_1q_2t_1t_2)_\infty}{(vt_1t_2)_\infty(vq_1q_2)_\infty}.
\end{eqnarray*}

Therefore,
\begin{align*}
T_1 +T_2 +T_3 &= \frac1{1 -t_1t_2}\cdot
 \frac{(v)_\infty
 (vq_1q_2t_1t_2)_\infty}{(vt_1t_2)_\infty(vq_1q_2)_\infty} \cdot
 \\
& \cdot \left[
 \frac{q_1q_2t_1t_2 -1}{(1-q_1t_1)(1-q_2t_2)}
 + \frac1{1-q_1t_1}
  \;{}_3 \Phi_2 \left ( \begin{array}{ccc}
  v, q_1t_1, vq_1q_2 \\
  vq_1, vq_1q_2t_1t_2
 \end{array}; v; t_2
 \right) \right. \\
 &\quad + \left. \frac1{1-q_2t_2}
  \;{}_3 \Phi_2 \left ( \begin{array}{ccc}
  v, q_2t_2, vq_1q_2 \\
  vq_2, vq_1q_2t_1t_2
 \end{array}; v; t_1
 \right) \right].
\end{align*}
Recall by  Lemma \ref{fun=value} that $\widehat{F}(q_1, t_1; q_2,
t_2) = (v)_\infty^{-1} \cdot \left \langle \widehat{B}_\la(q_1,
t_1) \widehat{B}_\la(q_2, t_2) \right \rangle_v.$ This together
with Lemma \ref{threeparts} proves Theorem~\ref{th:2point}.

In the case when $q_1q_2t_1t_2=1$, the above expression for $T_1
+T_2 +T_3$ can be further simplified to be
\begin{align} \label{hypersum}
 & \frac1{1 -t_1t_2} \cdot
 \frac{(v)_\infty^2}{(vt_1t_2)_\infty(vq_1q_2)_\infty} \cdot
 \nonumber
 \\
 &\cdot\left[
 \frac1{1-q_1t_1}
  \;{}_3 \Phi_2 \left ( \begin{array}{ccc}
  v, q_1t_1, vq_1q_2 \\
  vq_1, v
 \end{array}; v; t_2
 \right)
 +  \frac1{1-q_2t_2}
  \;{}_3 \Phi_2 \left ( \begin{array}{ccc}
  v, q_2t_2, vq_1q_2 \\
  vq_2, v
 \end{array}; v; t_1
 \right) \right]
 \nonumber \\
 %%%%%%
 %%%%%%
 &= \frac1{1 -t_1t_2} \cdot
 \frac{(v)_\infty^2}{(vt_1t_2)_\infty(vq_1q_2)_\infty} \cdot
 \nonumber
 \\
 & \cdot \left[
 \frac1{1-q_1t_1}
  \;{}_2 \Phi_1 \left ( \begin{array}{cc}
  q_1t_1, vq_1q_2 \\
  vq_1
 \end{array}; v; t_2
 \right)
 +  \frac1{1-q_2t_2}
  \;{}_2 \Phi_1 \left ( \begin{array}{cc}
   q_2t_2, vq_1q_2 \\
  vq_2
 \end{array}; v; t_1
 \right) \right]. \nonumber \\
\end{align}
Thanks to $q_1q_2t_1t_2 =1$, the two $(2,1)$-basic hypergeometric
series are of the form
$\;{}_2 \Phi_1 \left ( \begin{array}{cc}
   a,b \\
   c
 \end{array}; v; c/ab
 \right).$ Now by Heine's formula (cf. \cite[Appendix II.8]{GR})
 $$\;{}_2 \Phi_1 \left ( \begin{array}{cc}
   a,b \\
   c
 \end{array}; v; c/ab
 \right)
 =\frac{(c/a)_\infty (c/b)_\infty}{(c)_\infty (c/ab)_\infty},\quad |b|<1,|\frac{c}{ab}|<1,$$
 the expression (\ref{hypersum}) for $T_1
+T_2 +T_3$ becomes
\begin{align*}
 &\frac1{1 -t_1t_2} \cdot
 \frac{(v)^2_\infty}{(vt_1t_2)_\infty (vq_1q_2)_\infty}\cdot \\
 & \cdot \left[ \frac1{1 -q_1t_1}
 \frac{(vt_1^{-1})_\infty (q_2^{-1})_\infty}{(vq_1)_\infty
 (t_2)_\infty}
 +  \frac1{1 -q_2t_2}
 \frac{(vt_2^{-1})_\infty (q_1^{-1})_\infty}{(vq_2)_\infty
 (t_1)_\infty}
 \right].
\end{align*}
This together with Lemmas~\ref{fun=value} and \ref{threeparts}
completes the proof of Theorem~\ref{th:2special}.
\end{proof}
\begin{remark}
It follows from the proof above that the convergence of the
$2$-point function is guaranteed by assuming that $|t_1| <1$,
$|t_2|<1$,  $|vq_1q_2|<1$, $|v| <1$, and by excluding the values
for $t_i, q_i$ which make the denominators of the $(3,2)$-basic
hypergeometric series and other denominators in the above theorems
vanish.
\end{remark}
\section{A generalization via vertex operators}
\label{sec:general}
\subsection{$1$-point function of the zero-mode of a vertex operator}

Consider the Heisenberg algebra generated by $\text{I}$ and
$\mathfrak a_n, n\in\Z$ with the commutation relation (where
$\kappa$ is a constant):
\begin{eqnarray*}
  [ \mathfrak a_m, \mathfrak a_n] = \kappa m \delta_{m,-n} \text{I}.
\end{eqnarray*}
The Fock space $B$ is the irreducible representation of the
Heisenberg algebra generated by a (highest weight) vector $\vac$
such that $\text{I} \vac =\vac$ and $\mathfrak a_n \vac =0$ for $n
\ge 0$. The Fock space $B$ has a linear basis $\mathfrak a_{-\la}
:= \mathfrak a_{-\la_1}  \mathfrak a_{-\la_2} \cdots \vac$, where
$\la =(\la_1, \la_2, \ldots)$ runs over all partitions. Below we
identify $B$ with the ring of symmetric function $\Lambda$ by
identifying $\mathfrak a_{-\la}$ with the power-sum symmetric
functions $p_{\la}$.
%$B_{q,t} := B\otimes \mathbb Q(q,t)$ with $\Lambda_{q,t}$.

Introduce the following deformed vertex operator
\begin{eqnarray*}
 V(z;q_1, t_1, q_2, t_2) =
 \exp \left ( \sum_{k \ge 1} (q_1^k -q_2^k) \mathfrak a_{-k} \frac{z^k}k
  \right )
 \exp \left ( \sum_{k \ge 1} (t_2^k -t_1^k) \mathfrak a_{k} \frac{z^{-k}}k
  \right).
\end{eqnarray*}
Write
 $$V(z;q_1, t_1, q_2, t_2) =\sum_{m \in \Z}
V_m(q_1,q_2,t_1,t_2) z^m.$$

\begin{remark}
When $\kappa =1$, $q_2=t_2 =1$, and write $q=q_1$ and $t=t_1$, the
operator $V_0$ provides a vertex operator realization for
$\widehat{\mathfrak B}_{q,t}$:
\begin{eqnarray} \label{vo}
\widehat{\mathfrak B}_{q,t}
 =\frac{1}{(1-q)(1-t)} \cdot V_0(q,1,t,1).
 \end{eqnarray}
This formula in a $\la$-ring form (in different notations) appears
in the study of Macdonald polynomials by Garsia and Haiman
\cite[(73)]{GH}. In this sense, Theorem~\ref{th:1pointgeneral}
below is a generalization of Theorem~\ref{th:onepoint} (with
different proofs). The formula (\ref{vo}) for $t=q^{-1}$ is
equivalent to a formula of Lascoux and Thibon \cite[Prop.
3.3]{LT}.
\end{remark}
\begin{theorem}  \label{th:1pointgeneral}
 We have
 \begin{eqnarray*}
  \big \langle V_0(q_1,q_2,t_1,t_2) \big\rangle_v
  =   \left [
  \frac{(q_1 t_1 v)_\infty (q_2 t_2v)_\infty}{(q_1 t_2v)_\infty
  (q_2 t_1 v)_\infty}
  \right]^{\kappa}.
  \end{eqnarray*}
\end{theorem}

\begin{proof}
Let us denote $\Delta :=V_0(q_1,q_2,t_1,t_2)$. For a partition
$\la =(r^{m_r})_{r \ge 1}$ with $m_r$ parts equal to $r$,
$p_\la =  \prod_{r \ge 1}\mathfrak a_{-r}^{m_r} \vac.$ To compute
the trace $ \Tr_B v^{L_0} \Delta$, we will compute the projection
of $\Delta p_\la$ to the one-dimensional subspace $\C p_\la$ (with
respect to the basis $p_\mu$'s). A similar method has been also
used in \cite{DM}.
\begin{eqnarray*}
 && \text{projection of }\Delta p_\la\\
 &=& \sum_{\stackrel{(n_r)}{n_r \le m_r  \text{ for all } r}}
 \left( \prod_{r \ge 1} \frac{(q_1^r -q_2^r)^{n_r} \mathfrak a_{-r}^{n_r}}{r^{n_r}n_r!} \right)
 \left( \prod_{r \ge 1} \frac{(t_2^r -t_1^r)^{n_r} \mathfrak a_{r}^{n_r}}{r^{n_r}n_r!} \right)
 \cdot\prod_{r \ge 1} \mathfrak a_{-r}^{m_r}   \vac \\
 &=&  \sum_{\stackrel{(n_r)}{n_r \le m_r  \text{ for all } r}}
  \prod_{r \ge 1}  \frac{{m_r \choose n_r} r^{n_r}n_r! \kappa^{n_r}
  (q_1^r -q_2^r)^{n_r} (t_2^r -t_1^r)^{n_r}}{(r^{n_r}n_r!)^2}
  \cdot\prod_{r \ge 1} \mathfrak a_{-r}^{m_r}   \vac.
 \end{eqnarray*}
Therefore, we have
\begin{eqnarray*}
  \Tr_B (v^{L_0} \Delta )
  &=& \sum_{\stackrel{(m_r), (n_r)}{n_r \le m_r  \text{ for all } r}}
  \prod_{r \ge 1}  \frac{{m_r \choose n_r} r^{n_r}n_r! \kappa^{n_r}
  (q_1^r -q_2^r)^{n_r} (t_2^r -t_1^r)^{n_r} v^{rm_r}}{(r^{n_r}n_r!)^2} \\
  &=& \prod_{r \ge 1}
  \sum_{(n_r)} \frac{\left(\kappa
  (q_1^r -q_2^r)  (t_2^r -t_1^r) \right)^{n_r}}{r^{n_r}n_r!}
  \sum_{\stackrel{(m_r)}{m_r \ge n_r  \text{ for all } r}}
    {{m_r \choose n_r}   v^{rm_r}}.
\end{eqnarray*}
Using the simple binomial identity for $ n \ge 0$,
\begin{eqnarray*}
 \sum_{m \ge n} {m \choose n} x^m = \frac{x^n}{(1-x)^{1+n}},
\end{eqnarray*}
we have
\begin{eqnarray*}
  \Tr_B (v^{L_0} \Delta)
 &=& \prod_{r \ge 1} \sum_{(n_r)} \frac{\left(\kappa
  (q_1^r -q_2^r)  (t_2^r -t_1^r) \right)^{n_r}}{r^{n_r}n_r!}
 \frac{(v^{r})^{n_r}}{(1-v^r)^{1+n_r}} \\
 &=& (v )_\infty^{-1} \cdot \exp \left( \sum_{r \ge 1} \frac{
 \kappa v^r (q_1^r -q_2^r)  (t_2^r -t_1^r)}{r(1-v^r)}
 \right).
 \end{eqnarray*}
 Hence,
 \begin{eqnarray*}
 \langle \Delta \rangle_v &=& (v )_\infty \Tr_B (v^{L_0} \Delta) \\
 % &=& \exp \left( \sum_{r \ge 1} \frac{
 %\kappa v^r (q_1^r -q_2^r)  (t_2^r -t_1^r)}{1-v^r} \right) \\
 % %
  &=& \exp \left(  \sum_{r \ge 1} \sum_{n \ge 1} \frac{\kappa }{r}[
  (q_1t_2 v^n)^r +(q_2t_1 v^n)^r - (q_1t_1 v^n)^r -(q_2 t_2 v^n)^r
 ] \right) \\
 & =& \exp \left(\kappa  \sum_{n \ge 1}
 \big(\ln \,(1-q_1t_1 v^n) (1-q_2t_2 v^n)
  - \ln \, (1-q_1t_2 v^n)(1-q_2t_1 v^n) \big) \right) \\
 & =& \left [
  \frac{(q_1 t_1 v)_\infty (q_2 t_2v)_\infty}{(q_1 t_2v)_\infty
  (q_2 t_1 v)_\infty}
  \right]^{\kappa}.
 \end{eqnarray*}
\end{proof}

\subsection{The $n$-point function of a vertex operator}

It turns out that it is fairly easy to compute the $n$-point
function of the full vertex operator $V(z;s,t,u,w)$ in contrast to
the $n$-point function of its zero-mode (for $n \ge 2$). We first
recall a standard lemma.
\begin{lemma} \label{lem:single}
We have
\begin{align*}
 {\exp} & \left( \frac{(t_i^k-s_i^k) z_i^{-k}\mathfrak a_k}{k} \right)
 {\exp}\left ( \frac{(u_j^k-w_j^k) z_j^k \mathfrak a_{-k}}{k} \right) = \\
 &{\exp}\left(
\frac{\kappa (t_i^k-s_i^k)(u_j^k-w_j^k)z_j^kz_i^{-k}}{k}\right)
\times \\
 &\quad {\exp}\left ( \frac{(u_j^k-w_j^k) z_j^k \mathfrak a_{-k}}{k} \right)
 {\exp} \left( \frac{(t_i^k-s_i^k) z_i^{-k}\mathfrak a_k}{k} \right).
\end{align*}
\end{lemma}

\begin{lemma}  \label{lem:multi}
We have
\begin{align*}
\prod_{i=1}^n & V(z_i; s_i,t_i,u_i,w_i)
 = \\
& {\rm exp} \left(\kappa \sum_{k\ge 1}\frac{\sum_{1\le i<j\le
n}(t_i^k-s_i^k)(u_j^k-w_j^k)z_j^kz_i^{-k}}{k}\right)\times\\
&{\rm exp}\left( \sum_{k\ge 1}
\frac{ \sum_{i=1}^n(u_i^k-w_i^k)z_i^k \mathfrak a_{-k}}{k}\right)
{\rm exp}\left( \sum_{k\ge 1}
\frac{ \sum_{i=1}^n(t_i^k-s_i^k)z_i^{-k} \mathfrak a_k}{k}\right).
\end{align*}
\end{lemma}

\begin{proof}
Follows from applying Lemma~\ref{lem:single} repeatedly.
\end{proof}

\begin{theorem} \label{th:npoint}
 We have
  \begin{align}
 \left \langle  \prod_{i=1}^n V(z_i; s_i,t_i,u_i,w_i)
  \right \rangle_v  =&
    \prod_{1\le i<j\le n}
   \left [
   \frac{ (1 -t_i w_j z_i^{-1}z_j) (1 -s_i u_j z_i^{-1}z_j)}{ (1 -t_i u_j
   z_i^{-1}z_j) (1 -s_i w_j z_i^{-1}z_j)}
  \right]^{\kappa} \times \label{eq:prod1}\\
   &\quad \prod_{i,j=1}^{n}
   \left [
  \frac{(t_i w_j z_i^{-1}z_j)_\infty (s_i u_j z_i^{-1}z_j)_\infty}{(t_i u_j
   z_i^{-1}z_j)_\infty
  (s_i w_j z_i^{-1}z_j)_\infty}
  \right]^{\kappa}. \label{eq:prod2}
\end{align}
\end{theorem}

\begin{proof}
Using Lemma~\ref{lem:multi} and the projection technique as used
in the proof of Theorem~\ref{th:1pointgeneral}, the trace
%\begin{align*}
${\rm Tr} \left ( v^{L_0} \prod_{i=1}^n V(z_i; s_i,t_i,u_i,w_i)
\right)$
%\end{align*}
can be shown to be
\begin{align} \label{eq:inter}
(v )_\infty^{-1} \cdot {\exp} & \left( \kappa \sum_{k\ge 1}
\frac{\sum_{1\le i<j\le
n}(t_i^k-s_i^k)(u_j^k-w_j^k)z_j^kz_i^{-k}}{k}\right)\times \nonumber\\
{\exp}&\left(\kappa \sum_{k\ge 1}\frac{  v^k
\sum_{j=1}^n(u_j^k-w_j^k)z_j^k\cdot
\sum_{i=1}^n(t_i^k-s_i^k)z_i^{-k} }{k(1-v^k)}\right).
\end{align}

It is a simple algebraic manipulation to rewrite the first
exponential in (\ref{eq:inter}) as the product (\ref{eq:prod1})
 and the second exponential in (\ref{eq:inter}) as the product (\ref{eq:prod2}).
\end{proof}
\begin{remark}
Theorem~\ref{th:npoint} can be regarded as a generalization of
\cite[Theorem 3.1]{Mi}. It specializes when $n=1$ to
Theorem~\ref{th:1pointgeneral}.  For $n \ge 2$, the $n$-point
correlation function for a vertex operator differs from that for
the zero-mode of a vertex operator. While the correlation
functions for the zero-mode of a vertex operator has more direct
connections with other fields, it is much more difficult to
calculate.
\end{remark}

\section{Discussions}
\label{sec:discuss}

In this Note, we have formulated the $n$-point correlation
functions which are generalizations of \cite{BO}, and found closed
formulas when $n =1,2$. We then formulated and computed some
related $n$-point functions of vertex operators. In a way, this
Note raises more questions than we could answer. Let us list some
open problems and connections below:

\begin{enumerate}
\item The symmetric functions which are the eigenvectors for
$V_0(q_1,q_2,t_1,t_2)$ are common generalizations of the Macdonald
polynomials and Jack polynomials (with Jack parameter $\kappa$).
It is interesting to study them in detail and in particular to see
if they have Schur-positivity etc.

\item Calculate the $n$-point correlation functions for general
$n$. The simple closed formulas obtained in this Note for $n =1,2$
suggests a nice general solution, which will be a generalization
of the remarkable formula found in \cite{BO}.

\item The $n$-point functions of \cite{BO} afford geometric
interpretations in terms of Gromov-Witten theory of an elliptic
curve and Hilbert schemes of points on the affine plane. We
speculate that our $n$-point functions have similar
interpretations using equivariant $K$-theory formulations.

\item The function $\Bqt$ (after normalization) can be regarded as
a probability measure on the set of partitions, which generalizes
those studied actively in literature (cf. \cite{Ok} and the
references therein).
\end{enumerate}

{\bf Acknowledgment.} We thank George Andrews for a helpful email
correspondence. S.-J. C. is partially supported by an NSC grant of
the R.O.C., and he thanks University of Virginia for hospitality,
where part of this work was carried out. W.W. is partially
supported by NSA and NSF.

\end{document}